# Numerical Study of Nonlinear sine-Gordon Equation by using the Modified Cubic B-Spline Differential Quadrature Method

*H. S. Shukla*[1], *Mohammad Tamsir*[1*], *Vineet K. Srivastava*[2]

[1]Department of Mathematics & Statistics, DDU Gorakhpur University, Gorakhpur-273009, India
[2]ISRO Telemetry, Tracking and Command Network (ISTRAC), Bangalore-560058, India

**ABSTRACT**

In this article, we study the numerical solution of the one dimensional nonlinear sine-Gordon by using the modified cubic B-spline differential quadrature method (MCB-DQM). The scheme is a combination of a modified cubic B-spline basis function and the differential quadrature method. The modified cubic B-spline is used as a basis function in the differential quadrature method to compute the weighting coefficients. Thus, the sine-Gordon equation is converted into a system of ordinary differential equations (ODEs). The resulting system of ODEs is solved by an optimal five stage and fourth-order strong stability preserving Runge–Kutta scheme (SSP-RK54). The accuracy and efficiency of the scheme are successfully described by considering the three numerical examples of the nonlinear sine-Gordon equation having the exact solutions.

## 1. Introduction

The simulation of the physical phenomena arising in the various fields of science and engineering may be described by the partial differential equations.The partial differential equations that describe the nonlinear waves represented by the sine-Gordon equation are of the great importance in research, and the numerical methods that are suitable for the success of the simulation of waves has led to a constant demand for achieving the higher accuracy. The nonlinear one-dimensional sine–Gordon equation appears in the propagation of fluxons in Josephson junctions, dislocations in crystals, solid state physics, nonlinear optics, stability of fluid motions and the motion of a rigid pendulum attached to a stretched wire [1-3].

Consider the one-dimensional time dependent nonlinear sine-Gordon equation:

$$\frac{\partial^2 u(x,t)}{\partial t^2}=\frac{\partial^2 u(x,t)}{\partial x^2}-\sin(u(x,t)),\quad x\in[a,b],\quad t\geq 0, \tag{1.1}$$

with the initial conditions:

$$u(x,0) = f_1(x), \quad x \in [a,b],$$
$$\frac{\partial u}{\partial t}(x,0) = f_2(x), \quad x \in [a,b], \tag{1.2}$$
and Dirichlet boundary conditions:
$$u(a,t) = g_1(t), \quad u(b,t) = g_2(t), \quad t \geq 0. \tag{1.3}$$

In the recent years, various numerical approaches have been developed for the solution of Eq. (1.1) with the initial conditions (1.2) and boundary conditions (1.3). Guo et al. [4] proposed two different difference schemes, namely, explicit and implicit finite difference schemes with the advantage of conserved discrete energy and proved their stability and convergence. Fei & Vazquez [5] proposed two explicit conservative numerical schemes with their stability and convergence and showed that the new schemes are very accurate and fast. Bratsos & Twizell [7] have used a family of parametric finite difference methods to transform the nonlinear hyperbolic sine-Gordon equation into a system linear algebraic equation, while in [6] they used the method of lines to transform the nonlinear hyperbolic sine-Gordon equation into a first order, non-linear, initial value problem. Further, Bratsos [8] used the linearizing techniques on the nonlinear term of the sine-Gordon equation. Bratsos [9] also proposed a numerical scheme having second order accuracy in space and third-order accuracy in time by using a predictor-corrector method. He [10] further proposed a fourth order scheme for the solution of one dimensional sine-Gordon equation having second-order accuracy in space and fourth-order accuracy in time. Ablowitz et al. [11] investigated the numerical behavior of a double-discrete, completely integrable discretization of the Sine-Gordon equation and clarified the nature of the instability via an analytical investigation by the numerical experiments and further they [12] investigated how the order of accuracy and the symplectic property of the numerical scheme depends on the preservation of the integrable structure. Kaya [13] used a decomposition method to get the explicit and numerical solutions of the sine-Gordon equation in term of the convergent power series. It has been observed that this method performs extremely well in terms of accuracy, efficiency, simplicity, stability and reliability. Dehghan & Shokri [14] proposed a numerical scheme based on collocation points and thin plate splines radial basis function, while Dehghan & Mirzaei [15] used the boundary integral equation and dual reciprocity approach for the solution of one-dimensional sine-Gordon equation. The nonlinearity of the problem was dealt by the predictor–corrector scheme. The conservation of energy was also investigated by them [15]. Mohebbi & Dehghan [16] proposed a high order and accurate method based on the compact finite difference scheme for spatial components and the diagonally implicit Runge-Kutta-Nystrom (DIRKN) for the temporal components, which is fourth order accurate in both space, and time variables, and is unconditionally stable. Li-Min & Zong-Min [17] proposed multi-quadric quasi-interpolation approach; while a mesh free high accurate multi-quadric quasi-interpolation approach is proposed by Jiang & Wang [18]. Further, Jiang & Wang [19] used the spline function approximation and developed two implicit finite difference schemes for the numerical solution. Recently, Mittal & Bhatia [20] used a modified cubic B-spline collocation method. Differential quadrature technique can be seen in [21-26] for the solution of the linear and nonlinear one and two dimensional differential equations, while the cubic B-spline and sinc differential quadrature methods can be seen in references [27, 28]. Recently, Arora & Singh [29] proposed a new method, namely, "a modified cubic B-spline differential quadrature method (MCB-DQM) to solve the one dimensional Burgers' equation and checked its efficiency and accuracy, and observed that the implementation of this method is very easy, powerful and efficient as compared to other existing numerical methods for the Burgers' equation. Authors of [30-33] developed an optimal strong stability preserving (SSP) high order time discretization schemes for the numerical solution of the partial differential equations. SSP methods preserve the strong stability properties in any norm, semi norm or convex functional of the spatial discretization coupled with the first order Euler time stepping. They also discussed the description of the optimal explicit and implicit SSP Runge-Kutta and multistep methods.

In this work, we study the numerical solution of the one-dimensional nonlinear sine-Gordon equation using the MCB-DQM scheme. The efficacy and accuracy of the method is confirmed by taking the three examples of the Sine-Gordon equation having the exact solutions, which shows that the MCB-DQM results are acceptable and in good agreement with earlier studies available in the literature. Rest of the article is organized as follows: In Section 2, the modified cubic B-spline differential quadrature method is described. In Section 3, the implementation procedure is illustrated; in Section 4, three numerical examples are given in order to establish the applicability and accuracy of the method, while the Section 5 concludes our study.

## 2. Modified cubic B-spline differential quadrature method (MCB-DQM)

The differential quadrature method (DQM) was introduced by Bellman et al. [21] in 1972. This method approximates the spatial derivatives of a function using the weighted sum of the functional values at the certain discrete points. In DQM, the weighting coefficients are computed by using several kinds of test functions such as spline function [27], sinc function [28], Lagrange interpolation polynomials and Legendre polynomials [22-26] etc. This section re-describes the MCB-DQM [29] to complete our problem. It is assumed that the $N$ grid points $a = x_1 < x_2, ..... < x_N = b$ are distributed uniformly with $x_{i+1} - x_i = h$. The first and second order spatial derivatives [27] of $u(x,t)$ at $x_i$ are defined as:

$$\frac{\partial u(x_i,t)}{\partial x} = \sum_{j=1}^{N} w_{ij}^{(1)} u(x_j,t), \qquad i = 1,2,...,N \tag{2.1}$$

$$\frac{\partial^2 u(x_i,t)}{\partial x^2} = \sum_{j=1}^{N} w_{ij}^{(2)} u(x_j,t), \qquad i = 1,2,...,N, \tag{2.2}$$

where $w_{ij}^{(r)}$, $r = 1,2$ are the weighting coefficients of the $r$th-order spatial derivative.

The cubic B-spline basis functions [27] at the knots are defined as:

$$\varphi_j(x) = \frac{1}{h^3}\begin{cases} (x - x_{j-2})^3, & x \in (x_{j-2}, x_{j-1}) \\ (x - x_{j-2})^3 - 4(x - x_{j-1})^3, & x \in (x_{j-1}, x_j) \\ (x_{j+2} - x)^3, & x \in (x_{j+1}, x_{j+2}) \\ 0, & \text{else}, \end{cases} \tag{2.3}$$

where $\{\varphi_0, \varphi_1, ..., \varphi_N, \varphi_{N+1}\}$ is chosen in such a way that it forms a basis over $[a,b]$. The values of cubic B-splines and its derivatives at the nodal points are depicted in Table 1.

**Table 1:** Coefficients of the cubic B-spline $\varphi_j$ and its derivatives $\varphi_j', \varphi_j''$ at the node $x_j$.

| | $x_{j-2}$ | $x_{j-1}$ | $x_j$ | $x_{j+1}$ | $x_{j+2}$ |
|---|---|---|---|---|---|
| $\varphi_j(x)$ | 0 | 1 | 4 | 1 | 0 |
| $\varphi_j'(x)$ | 0 | $3/h$ | 0 | $-3/h$ | 0 |
| $\varphi_j''(x)$ | 0 | $6/h^2$ | $-12/h^2$ | $6/h^2$ | 0 |

The cubic B-spline basis functions are modified in such a way that the resulting matrix system of the equations is diagonally dominant. The modified cubic B-splines [29] are defined as:

$$\left.\begin{aligned} &\phi_1(x) = \varphi_1(x) + 2\varphi_0(x) \\ &\phi_2(x) = \varphi_2(x) - \varphi_0(x) \\ &\phi_j(x) = \varphi_j(x) \; for \; j = 3,...,N-2 \\ &\phi_1(x) = \varphi_{N-1}(x) - \varphi_{N+1}(x) \\ &\phi_1(x) = \varphi_N(x) + 2\varphi_{N+1}(x) \end{aligned}\right\}, \tag{2.4}$$

where $\{\phi_0, \phi_1, ..., \phi_N\}$ forms a basis over $[a,b]$.

The approximate value of the first order derivative at the $i^{th}$ node point is given by

$$\phi_k'(x_i) = \sum_{j=1}^{N} w_{ij}^1 \phi_k(x_j),\ k = 1,2,...,N. \tag{2.5}$$

From Eq. (2.3), (2.4) and Table 1, Eq. (2.5) is reduced into a tridiagonal system of equations:

$$A\,\vec{w}^{(1)}[i] = \vec{R}[i],\ for\ i = 1,2,...,N, \tag{2.6}$$

where $\vec{w}^{(1)}[i] = \left[w_{i1}^{(1)}, w_{i2}^{(1)}, ..., w_{iN}^{(1)}\right]^T$ is the weighting coefficient vector corresponding to $x_i$ .The coefficient matrix $A$ is given by

$$A = \begin{bmatrix} \phi_{1,1} & \phi_{1,2} & & & & & \\ \phi_{2,1} & \phi_{2,2} & \phi_{2,3} & & & & \\ & \phi_{3,2} & \phi_{3,3} & \phi_{3,4} & & & \\ & & \ddots & \ddots & \ddots & & \\ & & & \phi_{N-2,N-3} & \phi_{N-2,N-2} & \phi_{N-2,N-1} & \\ & & & & \phi_{N-1,N-2} & \phi_{N-1,N-1} & \phi_{N-1,N} \\ & & & & & \phi_{N,N-1} & \phi_{N,N} \end{bmatrix} \tag{2.7}$$

We notice that the coefficient matrix $A$, being diagonally dominant, is invertible. The tridiagonal system of linear equations (2.6) is solved by using Thomas algorithm. In this way, we obtained the weighting coefficients $w_{ij}^{(2)}, 1 \le i, j \le N$. The weighting coefficients $w_{ij}^{(2)}, 1 \le i, j \le N$, are determined by usingthe formula [16]:

$$\begin{cases} w_{ij}^{(r)} = r\left(w_{ij}^{(1)} w_{ii}^{(r-1)} - \dfrac{w_{ij}^{(r-1)}}{x_i - x_j}\right), for\ i \ne j \text{ and } i = 1,2,3,...,N; \quad r = 2,3,...,N-1 \\ w_{ii}^{(r)} = -\sum\limits_{j=1,\, j\ne i}^{N} w_{ij}^{(r)}, \text{ for } i = j. \end{cases} \tag{2.8}$$

## 3. Implementation of Method to Sine-Gordon equation

On substituting the approximate values of the spatial derivatives computed by using the MCB-DQM, Eq. (1.1) can be re-written as:

$$\frac{\partial^2 u(x_i,t)}{\partial t^2} = \sum_{j=1}^{N} w_{ij}^{(2)} u(x_j) - \sin(u_i), \qquad i = 1,2,...,N. \tag{3.1}$$

Eq. (3.1) is reduced into a set of first-order ODEs in time:

$$\frac{d^2u_i}{dt^2}=L\left(u_i\right),\ \ i=1,2,...,N, \tag{3.2}$$

subject to the initial conditions (1.2), where $L$ is a nonlinear differential operator. We use the well-known SSP-RK54 scheme [32] to solve the resulting system (3.2) and consequently the solution $u(x,t)$ is obtained at the required time level.

## 4. Numerical results and discussion

In section, we consider three numerical examples to provide the MCB-DQM numerical solutions of the nonlinear Eq. (1.1) having the exact solutions. The accuracy and efficiency of the MCB-DQM is measured by evaluating the $L_2$, $L_\infty$, and RMS error norms defined as:

$$\left.\begin{aligned} &L_2=\sqrt{h\sum_{j=0}^{N}\left|u_j^{exact}-u_j^{numerical}\right|^2};\ L_\infty=\max_{0\le j\le N}|u_j^{exact}-u_j^{numerical}| \\ &\text{and}\quad \text{RMS}=\frac{1}{N+1}\sqrt{\sum_{j=0}^{N}\left|u_j^{exact}-u_j^{numerical}\right|^2} \end{aligned}\right\}. \tag{4.1}$$

**Example 4.1**: Consider the sine-Gordon equation (1.1) with the initial conditions:

$$u\left(x,0\right)=0\ \text{ and }\ u_t(x,0)=4\sec h(x). \tag{4.2}$$

The exact solution [14, 16-20] of Eq. (1.1) with the initial conditions (4.2) is given by:

$$u(x,t)=4\tan^{-1}(\sec h(x)t). \tag{4.3}$$

The boundary conditions are extracted from the exact solutions (4.3).

Firstly we compute the numerical results in the domain [−1, 1] using the parameters: $\Delta t=0.0001$ and the space step size $h=0.04$ and compare our results with the results available in the literature [14, 20].**Table 2** shows the $L_2$ and $L_\infty$ error norms at the different time levels. It can be seen from **Table 2** that the MCB-DQM results are in good agreement with those of [14, 20]. Next, we compute the numerical results in the domain [−2, 2] with $\Delta t=0.01$ and $h=0.01$ as considered by Li-Min & Zong-Min [17], Jiang & Wang [18], and Mittal & Bhatia [20]. The $L_\infty$ and RMS norms at the different time levels are shown in **Table 2**. From **Table 3**, it can be observed that the MCB-DQM results are better than the results obtained in [17, 18, 20]. **Fig. 1** shows the comparison between the numerical and exact solutions at $t=1$ with $h=0.04$ and $\Delta t=0.0001$ while **Fig. 2** shows the physical behavior of the numerical solution in $3D$ with $h = \Delta t = 0.01$. and the contour form of the example 4.1, qualitatively.

**Table 2:** Comparison of the $L_2$ and $L_\infty$ error norms of the example 4.1 in the domain [-1, 1] with $h = 0.04$ and $\Delta t = 0.0001$ at the different time levels.

| $t$ | $L_2$ | $L_\infty$ | $L_2$ | $L_\infty$ | $L_2$ | $L_\infty$ |
|---|---|---|---|---|---|---|
| | MCB-DQM | | Mittal & Bhatia [20] | | Dehghan & Shokri [14] | |
| 0.25 | $2.43\times10^{-6}$ | $5.46\times10^{-6}$ | $1.18\times10^{-5}$ | $2.32\times10^{-5}$ | $3.91\times10^{-5}$ | $5.89\times10^{-6}$ |
| 0.50 | $5.54\times10^{-6}$ | $7.39\times10^{-6}$ | $4.19\times10^{-5}$ | $4.11\times10^{-5}$ | $1.30\times10^{-4}$ | $2.01\times10^{-5}$ |
| 0.75 | $6.45\times10^{-6}$ | $7.40\times10^{-6}$ | $7.78\times10^{-5}$ | $1.02\times10^{-4}$ | $2.35\times10^{-4}$ | $3.63\times10^{-5}$ |
| 1.0 | $7.84\times10^{-6}$ | $8.75\times10^{-6}$ | $1.30\times10^{-4}$ | $1.64\times10^{-4}$ | $3.27\times10^{-4}$ | $5.07\times10^{-5}$ |

**Table 3:** Comparison of the $L_\infty$ and RMS error norms of the example 4.1 in the domain [-2, 2] with $h = \Delta t = 0.01$ at the different time levels.

| $t$ | | MCB-DQM | Li-Min & Zong-Min [17] | Jiang & Wang [18] | Mittal & Bhatia [20] |
|---|---|---|---|---|---|
| 0.2 | $L_\infty$ | $1.46\times10^{-6}$ | $9.25\times10^{-5}$ | $2.50\times10^{-5}$ | $2.26\times10^{-5}$ |
| | RMS | $2.54\times10^{-8}$ | $1.76\times10^{-5}$ | $6.55\times10^{-7}$ | $2.69\times10^{-7}$ |
| 0.4 | $L_\infty$ | $2.97\times10^{-6}$ | $1.62\times10^{-4}$ | $4.20\times10^{-5}$ | $7.52\times10^{-5}$ |
| | RMS | $5.67\times10^{-8}$ | $1.62\times10^{-4}$ | $1.15\times10^{-6}$ | $1.19\times10^{-6}$ |
| 0.6 | $L_\infty$ | $4.32\times10^{-6}$ | $3.73\times10^{-4}$ | $6.54\times10^{-5}$ | $1.55\times10^{-4}$ |
| | RMS | $9.07\times10^{-8}$ | $1.65\times10^{-4}$ | $1.55\times10^{-6}$ | $2.96\times10^{-6}$ |
| 0.8 | $L_\infty$ | $5.46\times10^{-6}$ | $6.24\times10^{-4}$ | $4.01\times10^{-4}$ | $2.59\times10^{-4}$ |
| | RMS | $1.25\times10^{-7}$ | $2.98\times10^{-4}$ | $3.92\times10^{-6}$ | $5.72\times10^{-6}$ |
| 1.0 | $L_\infty$ | $6.33\times10^{-6}$ | $8.49\times10^{-4}$ | $1.53\times10^{-3}$ | $3.84\times10^{-4}$ |
| | RMS | $1.58\times10^{-7}$ | $4.37\times10^{-4}$ | $1.56\times10^{-5}$ | $9.56\times10^{-6}$ |

**Table 4:** Order of convergence for the example 4.2 in the domain [-3, 3] at $t = 1$.

| $h$ | $L_2$ | Order | $L_\infty$ | order |
|---|---|---|---|---|
| 0.04 | $3.331197\times10^{-5}$ | - | $3.743029\times10^{-5}$ | - |
| 0.02 | $8.480206\times10^{-6}$ | 1.974 | $9.640778\times10^{-6}$ | 1.957 |
| 0.01 | $2.087315\times10^{-6}$ | 2.023 | $2.385016\times10^{-6}$ | 2.015 |
| 0.005 | $4.660351\times10^{-7}$ | 2.163 | $5.340646\times10^{-7}$ | 2.159 |

**Example 4.2**: In this example, we consider the sine-Gordon equation (1.1) in the range $-3 \le x \le 3$ with the initial conditions:

$$\left.\begin{aligned} u(x,0) &= 4\tan^{-1}\left(\exp(\gamma x)\right) \\ u_t(x,0) &= \frac{-4c\gamma\exp(\gamma x)}{1+\exp(2\gamma x)} \end{aligned}\right\}, \tag{4.4}$$

where $c$ is the velocity of solitary wave and $\gamma = \dfrac{1}{\sqrt{1-c^2}}$.

The exact solution [16, 19-20]of Eq. (1.1) with the initial conditions (4.4) is given as:

$$u(x,t) = 4\tan^{-1}\left(\exp\left(\gamma(x-ct)\right)\right). \tag{4.5}$$

The boundary conditions are taken from the exact solutions (4.5).

The numerical results of the example 4.2 are computed in the computational domain $[-3, 3]$ with the parameters: $c = 0.5$, $\Delta t = 0.0001$ and $h = 0.04$ . Computed results are compared with the results obtained by Dehghan & Shokri [14], and Mittal & Bhatia [20]. **Table 5** shows the $L_2$ and $L_\infty$ error norms at the different time levels, quantatively. From **Table 5,** it can be seen that the MCB-DQM results are in good agreement with those of [14, 20]. In addition,**Table 4** shows the $L_2$ and $L_\infty$ norms and the order of convergence for $\Delta t = 0.01$ and for the different values of $h$ at $t = 1.0$, while **Table 6** shows the order of convergence for $c = 0.5$, $\Delta t = 0.0001$ at $t = 1.0$. It is clear that the rate of convergence of the described scheme is quadratic. **Fig. 3** shows the comparison between the numerical and exact solutions at $t = 1$ with $h = 0.08$ and $\Delta t = 0.0001$.The concentration profiles of the numerical solution in $3D$ and the contour form are depicted in **Fig. 4** for $0 \le t \le 1$ in the domain [-3, 3] with $h = 0.04$.

**Table 5:** Comparison of the $L_\infty$ and $L_\infty$ error norms of the example 4.2 in the domain [-3, 3] with $h = 0.04$, $\Delta t = 0.0001$ and $c = 0.5$ at the different time levels.

| $t$ | MCB-DQM | | Dehgam & Shokri [14] | | Mittal & Bhatia [20] | |
|---|---|---|---|---|---|---|
| | $L_2$ | $L_\infty$ | $L_2$ | $L_\infty$ | $L_2$ | $L_\infty$ |
| 0.25 | $5.67\times10^{-6}$ | $9.61\times10^{-6}$ | $1.76\times10^{-5}$ | $4.95\times10^{-6}$ | $3.66\times10^{-5}$ | $4.90\times10^{-5}$ |
| 0.50 | $8.39\times10^{-6}$ | $1.10\times10^{-5}$ | $4.31\times10^{-5}$ | $8.42\times10^{-6}$ | $9.00\times10^{-5}$ | $7.55\times10^{-5}$ |
| 0.75 | $1.05\times10^{-5}$ | $1.26\times10^{-5}$ | $8.25\times10^{-5}$ | $1.65\times10^{-5}$ | $1.60\times10^{-4}$ | $1.43\times10^{-4}$ |
| 1.0 | $1.24\times10^{-5}$ | $1.44\times10^{-5}$ | $1.27\times10^{-4}$ | $2.51\times10^{-5}$ | $2.27\times10^{-4}$ | $2.10\times10^{-4}$ |

**Table 6:** Order of convergence for the example 4.2 in the domain [-3, 3] at $t = 1$.

| $h$ | $L_2$ | Order | $L_\infty$ | order |
|---|---|---|---|---|
| 0.04 | $1.235453\times10^{-5}$ | - | $1.439969\times10^{-5}$ | - |
| 0.02 | $3.208207\times10^{-6}$ | 1.945 | $3.820306\times10^{-6}$ | 1.914 |
| 0.01 | $8.168641\times10^{-7}$ | 1.974 | $9.834778\times10^{-7}$ | 1.958 |
| 0.005 | $2.05831\times10^{-7}$ | 1.987 | $2.493047\times10^{-7}$ | 1.980 |

**Example 4.3**: In this example, we consider the Sine-Gordon equation (1.1) in the range $-10 \le x \le 10$ with the initial conditions:

$$\left.\begin{aligned} u(x,0) &= 0 \\ u_t(x,0) &= 4\overline{\gamma}\sec h(\overline{\gamma}x) \end{aligned}\right\}, \tag{4.6}$$

where $c$ is the velocity of solitary wave and $\overline{\gamma} = \dfrac{1}{\sqrt{1+c^2}}$.

The exact solution [10, 20] of Eq. (1.1) with the initial conditions (4.6) is given by:

$$u(x,t) = 4\tan^{-1}\left(c^{-1}\sin(\overline{\gamma}ct)\sec h(\overline{\gamma}x)\right). \tag{4.5}$$

The boundary conditions are extracted from the exact solutions (4.5).

The numerical results for the example 4.3 are computed in the domain [−10, 10] using the parameters: $c = 0.5$, $\Delta t = 0.001$ and $h = 0.01$ . Computed results are compared with the results obtained by Brastos [10], and Mittal & Bhatia [20]. **Table 7** show the $L_2$ and $L_\infty$ norms at the different time levels. From **Table 8,** it can be seen that the MCB-DQM results are in good agreement with those of [10, 20]. **Fig. 5** shows the comparison between the numerical and exact solutions at $t = 20$ with $h = 0.04$ and $\Delta t = 0.001$. The concentration profiles of the numerical solution in $3D$ and the contour form are shown in **Fig. 5** for $0 \le t \le 20$ in the domain [-10, 10].We notice that the concentration profile of the numerical solution is quite similar to Mittal & Bhatia [20].

**Table 7:** Comparison of the $L_\infty$ and $L_\infty$ error norms of the example 4.3 in the domain [-10, 10] with $h = 0.01$ , $\Delta t = 0.001$ and at the different time levels.

| $t$ | MCB-DQM | | Mittal & Bhatia [20] | | Brastos [10] |
|---|---|---|---|---|---|
| | $L_2$ | $L_\infty$ | $L_2$ | $L_\infty$ | $L_\infty$ |
| 1 | $1.866\times10^{-9}$ | $2.318\times10^{-9}$ | $2.564\times10^{-5}$ | $1.818\times10^{-5}$ | $1.276\times10^{-4}$ |
| 10 | $5.474\times10^{-9}$ | $5.234\times10^{-9}$ | $8.850\times10^{-5}$ | $5.228\times10^{-5}$ | $1.912\times10^{-4}$ |
| 20 | $9.800\times10^{-9}$ | $5.471\times10^{-9}$ | $1.713\times10^{-4}$ | $9.438\times10^{-5}$ | $2.519\times10^{-4}$ |

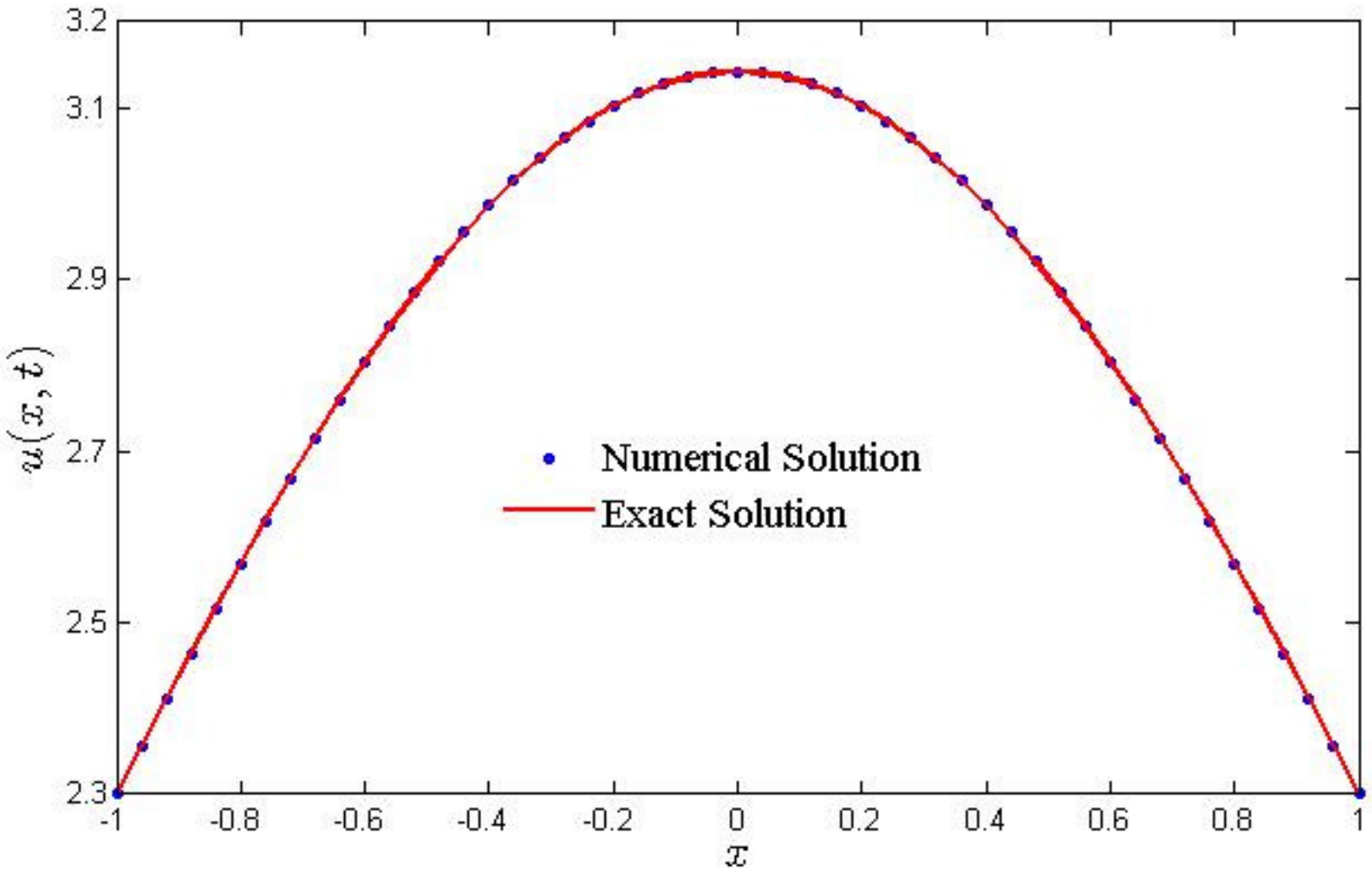

**Fig. 1:** Comparison between the numerical and exact solutions for the example 4.1 at $t = 1$ with $h = 0.04$ and $\Delta t = 0.0001$ in the range $-1 \le x \le 1$.

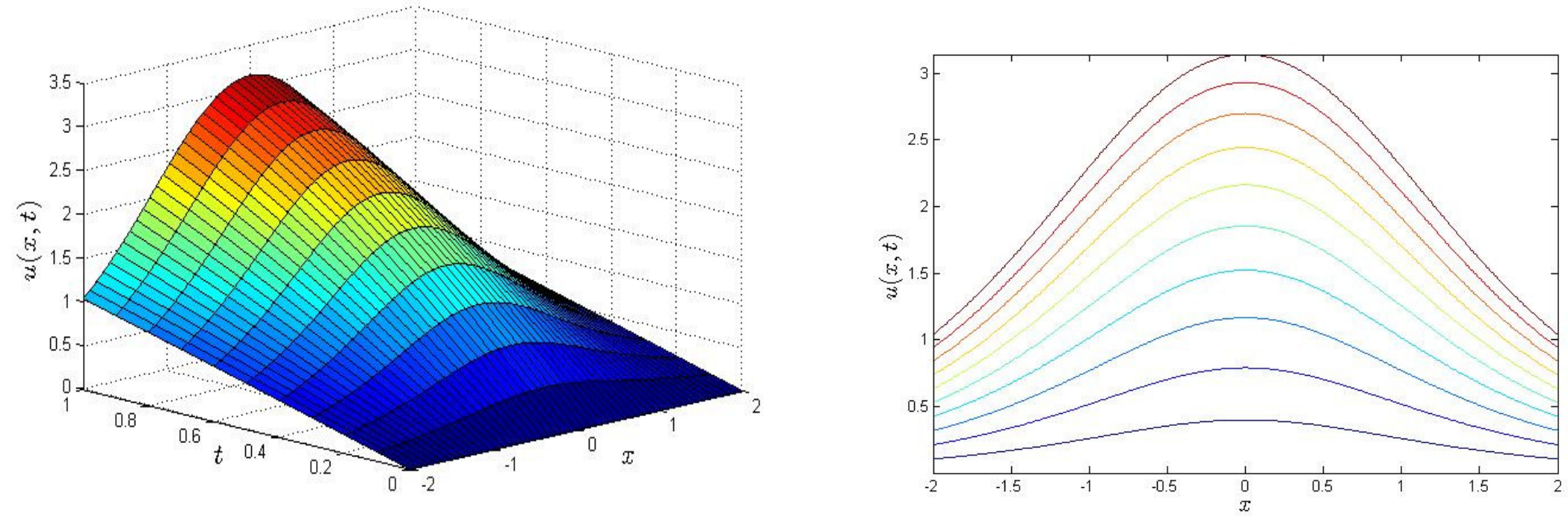

**Fig. 2:** Physical behavior of the numerical solution in $3D$ (left) and the contour form (right) of the example 4.1 for $0 \le t \le 1$ in the range $-2 \le x \le 2$.

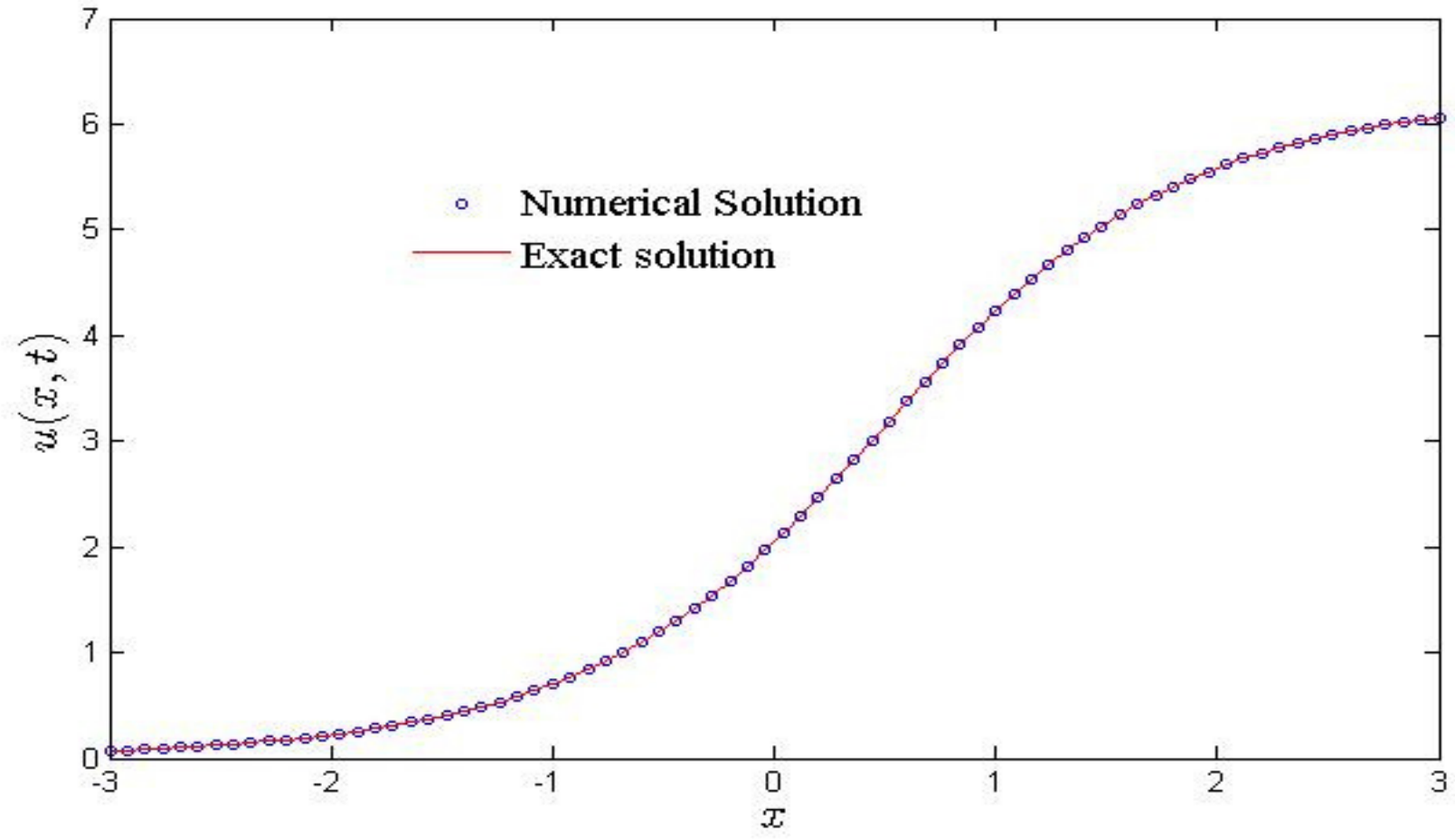

**Fig. 3:** Comparison between the numerical and exact solutions of the example 4.2 at $t = 1$ with $h = 0.08$ and $\Delta t = 0.0001$.

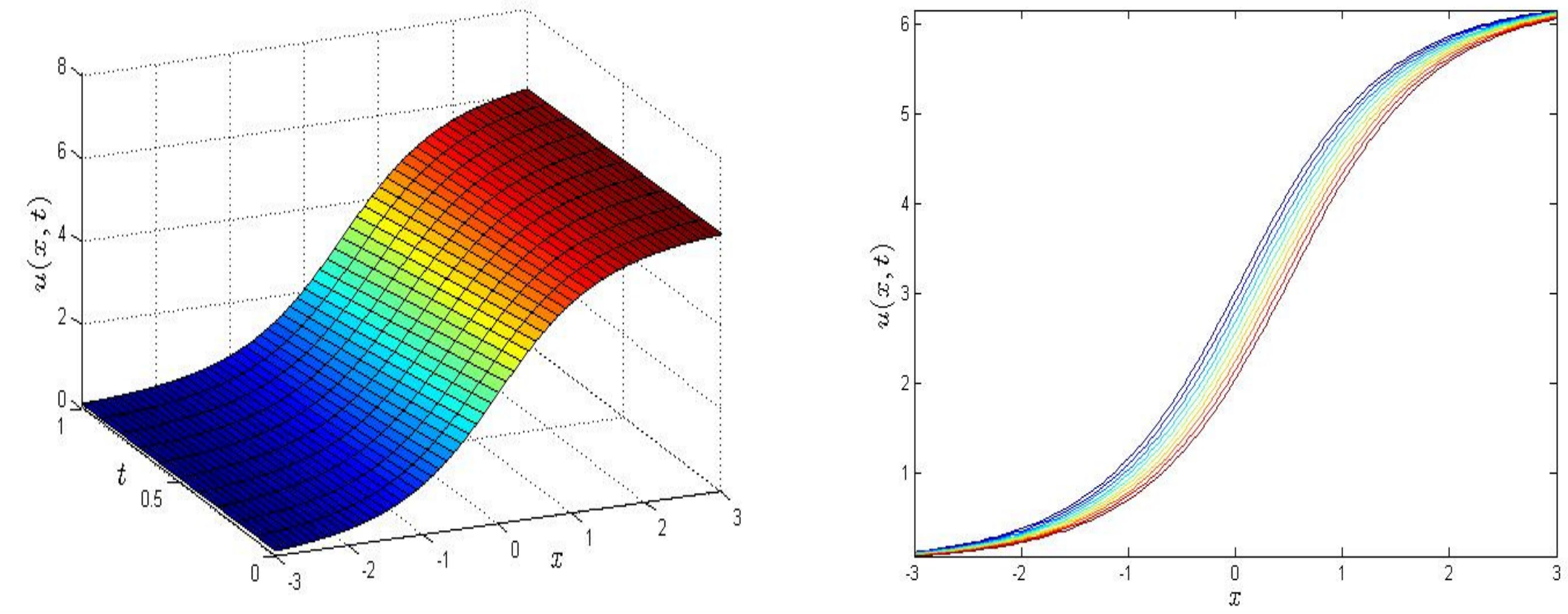

**Fig. 4:** Physical behavior of the numerical solution in $3D$ (left) and the contour form (right) of the example 4.2 for $0 \le t \le 1$ in the range $-3 \le x \le 3$.

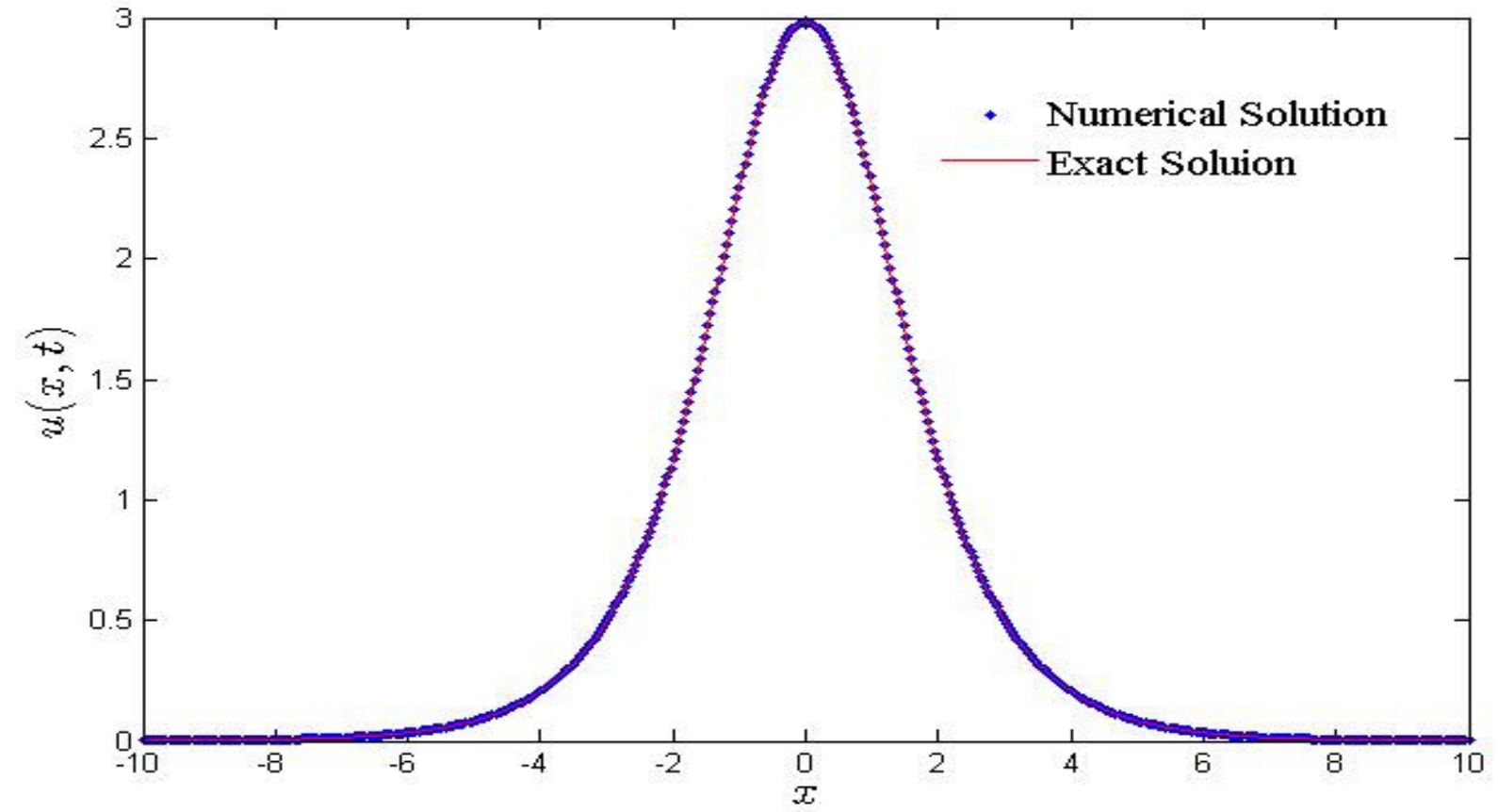

**Fig. 5:** Comparison between the numerical and exact solutions of the example 4.3 at $t = 20$ with $h = 0.04$ and $\Delta t = 0.001$.

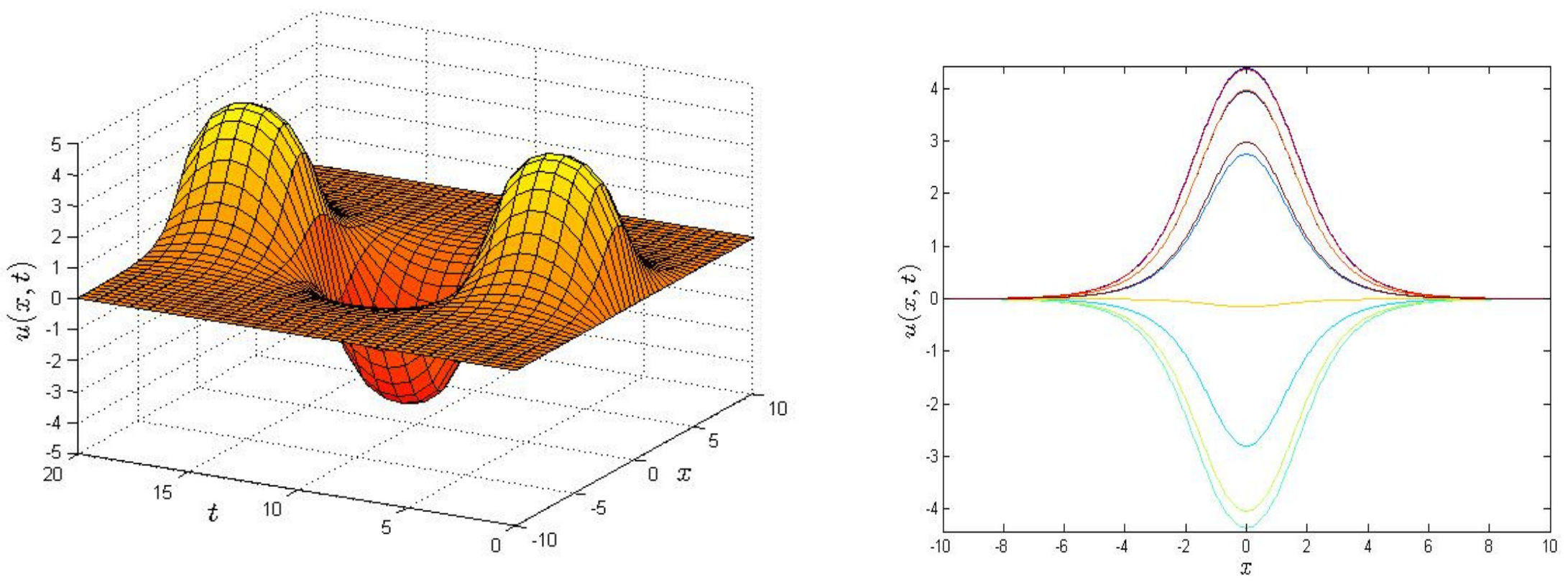

**Fig. 6:** Physical behavior of the numerical solution in 3 $D$ (left) and the contour form (right) of the example 4.3 for $0 \le t \le 20$ in the range $-10 \le x \le 10$.

## 5. Conclusions

This article explored the utility of a composite scheme: MCB-DQM, in space with SSP-RK54 scheme, in time for solving the one dimensional sine-Gordon equation. The efficiency, accuracy and reliability of the scheme are illustrated through the three numerical examples. It is evident that the error norms in all the examples are very small; an excellent numerical approximation to the exact solution is achieved. It is demonstrated that the MCB-DQM solution of the Sine-Gordon equation is very close to the exact solution and performs better than other numerical results available in the literature.